# Infinite-Dimensional System Signature: System Signature Under the Repairable Principle

Y. GÜRAL, M. GÜRCAN

*Abstract*— Investigation of the reliability of technical systems is one of the application areas of stochastic processes. The reliability of a technical system is based on two main elements. The first is the connection type of the system, and the second is the distribution of the working times of the components consisting of the system. In this study, system signatures and their reliability will be calculated under the repairable principle of parallel and serial systems consisting of two components. Although there are a limited number of studies in the literature for repairable systems, there is no study on creating the signature of repairable systems. In technical systems where there is no repair principle, although the system signature has limited components, the technical systems working under the repairable principle, have infinite components of the system signature. While creating the system signature, the probability that the working time of the component that is in the state of the system failure is greater than the repair time was defined as the parameter $\xi$. In the application part of the study, under the principle of repair, the system signature, and the reliability of the system were be successfully calculated.

## I. Introduction

One of the most important problems in engineering applications is the good design of the system design. This design is an important criterion that affects the performance and cost of the system. Generally, many technical systems consist of designing parallel and serially connected components in different ways. For example, the structure-function of k out of n system can be easily written in terms of parallel and serial connected systems.

Let $X_1, \cdots, X_n$ ve $C(n,k) = r$ be the components that constitute the technical system and let $A_1, \cdots, A_r$ be the k-elements subsets of this system. In this case, the structure-function of $k$ out of $n$ system can be written as follows,

$$min\{maxA_1, \cdots, maxA_r\}$$

where $A_j = \{X_{j_1}, \cdots, X_{j_k}\}$, $j_1, \cdots, j_k \in \{1, \cdots, n\}$ and $maxA_j = max\{X_{j_1}, \cdots, X_{j_k}\}$. Similarly, the structure-function of consecutive $k$ out of $n$ system can be shown as follows,

$$min\{maxB_1, \cdots, maxB_s\}$$

where $s = n - k + 1$, $B_j = \{X_j, X_{j+1}, \cdots, X_{j+k-1}\}$.

In calculating the reliability of technical systems, the number and location of the components that are working and failed are very important. For example, the reliability of parallel-connected systems is highest among all technical systems, and the reliability of serial connected technical systems is the lowest among all technical systems. Therefore, system signature is an important measure for the reliability of the system. The system signature is a probability vector, and the $k$-$th$ component of the system signature indicates the probability that the kth failed component will fail the system. The system signature was first used by Samaniego in the calculation of the reliability of the technical system as follows [1-3],

$$P(U > t) = \sum_{i=1}^{n} s_i Pr\{U_{i:n} > t\} \qquad (1)$$

where $S = (s_1, s_2, \ldots, s_n)$ is system signature, $U$ is the lifetime of the system and $U_{i:n}$ is the lifetime of the $i$-$th$ failed component.

Navarro et al. showed how the identification of the system signature can be obtained, which will then make it possible to compare consistent systems with a different number of components. They also stated that these results are also valid for mixed systems [4]. Kochar et al. proposed a new approach based on a standardizing assumption that sets the compared systems on an equal footing. They showed that this approach provides a framework from which very strong conclusions can be drawn regarding lifetime distributions of competing systems [5]. Eryilmaz presented a general framework to show the usefulness of the system signature in the reliability analysis of repairable systems with a general structure [7]. In any of these studies examined in the literature, it was not possible to create a system signature under the principle of repair. Therefore, the system signature is not used for the reliability of repairable systems. In this study, how a system signature of a repairable system can be created in infinite dimensions is presented. In this way, the system reliability of a repairable system can be calculated with the help of the system signature.

## II. System Signature and Reliability Under the Repairable Principle of the Parallel Connected Technical System

First, let's examine the system of $max\{X_1, X_2\}$, which consists of two components. For such a designed system to fail under the repairable principle, the working time of the failed component must be shorter than the repair time of the failed

Y. GÜRAL is with Firat University Faculty of Science Department of Statistics, Elazig, TURKEY. (e-mail: ygural@firat.edu.tr).

M. GÜRCAN is with Firat University Faculty of Science Department of Statistics, Elazig, TURKEY. (e-mail: mgurcan@firat.edu.tr).

component when it is in repair. Let $t_1, t_2, \cdots$ be the moments of failure of the components. Let $Y$ and $G(y)$ be the repair time and distribution function respectively. While $(j-1)$-$th$ component is in repair, the probability of working of the system $\xi_j = \Pr\{T_j > y\}$ is calculated as follows,

$$\xi_j = \int_0^\infty \Pr\{T_j > y, y \geq 0\} dy$$
$$= \int_0^\infty \Pr\{T_j > y, y \geq 0\} dy$$
$$= \int_0^\infty \Pr\{T_j > y\} dG(y)$$

where $T_1 = t_1$, $T_j = t_j - t_{j-1}$, $j = 2,3,\cdots$ and $Y$ denote the random variable repair time.

While $(j-1)$-$th$ component is in repair, the probability of system failure will be $1 - \xi_j$. Since the system is connected in parallel, the probability of the first failed component to fail the system is zero. The probability of the second failed component to fail the system is $1 - \xi_2$. For the third failed component to fail the system, the second failed component must not fail the system at the same time. Thus, the probability of the third failed component to fail the system is $(1 - \xi_3)\xi_2$. The probability of the $j$-$th$ failed component to fail the system is $(1 - \xi_j)\xi_2 \cdots \xi_{2+j-3}$, $j \geq 3$. Accordingly, the components of the system signature of the technical system are as follows,

$$s_1 = 0$$
$$s_2 = 1 - \xi_2$$
$$s_j = (1 - \xi_j)\xi_2 \cdots \xi_{2+j-3}, \ j \geq 3$$

*Remark 1:* Let's assume that $T_j \sim \text{Exponential}(\lambda)$, $Y \sim \text{Exponential}(\mu)$. Then $\xi_2 = \xi_3 = \cdots = \xi$ and the components of the system signature are obtained as follows,

$$s_1 = 0 \quad (2)$$
$$s_2 = 1 - \xi \quad (3)$$
$$s_j = (1 - \xi)\xi^{j-2}, \ j \geq 3 \quad (4)$$

where $\xi = \mu/(\lambda + \mu)$.

*Remark 2:* When $\xi_2 = \xi_3 = \cdots = \xi = 0$, the system becomes the non-repairable system. When $\xi_2 = \xi_3 = \cdots = \xi > 0$, the system becomes the repairable system.

*Remark 3:* The system signature of n-component and parallel-connected technical system under the repairable principle is obtained as follows,

$$s_1 = s_2 = \cdots = s_{n-1} = 0, s_n = 1 - \xi,$$
$$s_{n+1} = 1 - \xi, s_j = (1 - \xi)\xi^{j-n}, \ j \geq n$$

where failure times are independent and identically distributed (i.i.d).

### III. SYSTEM SIGNATURE AND RELIABILITY UNDER THE REPAIRABLE PRINCIPLE OF THE SERIAL CONNECTED TECHNICAL SYSTEM

Let's examine the system signature of the $min\{X_1, X_2\}$ system consisting of two components under the repairable principle. To examine the working principle of a serially connected technical system under the condition that it can be repaired, a replacement component must first be in the system.

Naturally, such a system consists of two series-connected components and a replacement component. In this case, the probability of the first component failing the system is zero. For the second failed component to fail the system, the repair time must be greater than the working times of the working components. This possibility is written as follows,

$$\Pr\{min(X_2, Y) > Z_1\} = \Pr\{X_2 > Z_1, Y > Z_1\}$$
$$= Pr\{X_2 > Z_1\}Pr\{Y > Z_1\}$$
$$= (\mu/(\lambda + \mu))^2$$

where $X_1$ and $Z_1$ be the first failed component and repair time, respectively. The working times and repair times of the components were accepted as i.i.d with $\lambda$ parameter exponential distribution and $\mu$ parameter exponential distribution, respectively.

Similarly, it can be obtained the probability of the system working during ongoing failure. In this case, the system signature $\xi = (\mu/(\lambda + \mu))^2$ of the serial connected system will be similar to the parallel-connected system.

*Remark 4:* The structure-function of the repairable technical system with serially connected two components can be written as follows,

$$max\{min(max(X_1, Y), X_2), min(X_1, max(X_2, Y))\}$$

where $Y$ is replacement component.

*Remark 5:* The system signature of the technical system with the structure-function $min(max(X_1, X_2), X_3)$ can be obtained as follows,

$$S = (q, p(1 - \xi), p(1 - \xi)\xi, p(1 - \xi)\xi^2, \cdots)$$

where $q$ is the probability of failure of the component connected in series to the technical system, and $p$ is the probability of working.

### IV. ORDER STATISTICS OF FAILED COMPONENTS UNDER REPAIRABLE PRINCIPLE

Let $T_k$ be the lifetime of the $k$-$th$ failed component under the repairable principle. In this case, the possible states for the component to be failed in the $k$-$th$ order are as follows.

The component may fail for the first time: $\Pr\{T_k > t\} = 1 - F(t), F(t) \sim Exp(\lambda)$

The component may fail for the second time: In this case, since the component has been repaired once before, under the condition that the repair time is $Exp(\mu)$,

$$\Pr\{T_k > t\} = Pr\{Erlang(2, \lambda) + Erlang(1, \mu) > t\} \quad (5)$$

Similarly, the component may have been failed $(k-1)$-$th$ times: In this case, it will be repaired $k - 1$ times,

$$\Pr\{T_k > t\} = Pr\{Erlang(k, \lambda) + Erlang(k - 1, \mu) > t\} \quad (6)$$

where $Erlang(k, \lambda)$ is meant the random variable with $k$ degrees of freedom Erlang distribution with parameter $\lambda$. The $k$ degree of freedom Erlang distribution is derived from the convolution of independent $k$ exponential distributions.

## V. ILLUSTRATIVE EXAMPLE

Now, for a 2-component parallel-connected repairable system, an application has been made for the case where the working time of the components is exponential with $\lambda$ parameter and the repair time is exponential with $\mu$ parameter.

The values of the system signature vector and the order statistics were obtained by using (2), (3), (4) and (6), respectively. With these obtained values, the reliability values in tables 1-3 were calculated using (1).

TABLE I
SYSTEM RELIABILITY FOR $\lambda = 0,1$ AND $\mu = 0,2$.

| t | Reliability | t | Reliability |
|---|---|---|---|
| 1 | 0,9501 | 21 | 0,5504 |
| 2 | 0,9279 | 22 | 0,5323 |
| 3 | 0,9067 | 23 | 0,5145 |
| 4 | 0,8861 | 24 | 0,4972 |
| 5 | 0,8659 | 25 | 0,4802 |
| 6 | 0,8459 | 26 | 0,4636 |
| 7 | 0,8259 | 27 | 0,4475 |
| 8 | 0,8058 | 28 | 0,4317 |
| 9 | 0,7857 | 29 | 0,4164 |
| 10 | 0,7656 | 30 | 0,4015 |
| 11 | 0,7454 | 31 | 0,3870 |
| 12 | 0,7253 | 32 | 0,3729 |
| 13 | 0,7052 | 33 | 0,3592 |
| 14 | 0,6851 | 34 | 0,3459 |
| 15 | 0,6652 | 35 | 0,3330 |
| 16 | 0,6455 | 36 | 0,3206 |
| 17 | 0,6259 | 37 | 0,3085 |
| 18 | 0,6066 | 38 | 0,2968 |
| 19 | 0,5875 | 39 | 0,2855 |
| 20 | 0,5688 | 40 | 0,2746 |

TABLE II
SYSTEM RELIABILITY FOR $\lambda = 0,0667$ AND $\mu = 0,1$.

| t | Reliability | t | Reliability |
|---|---|---|---|
| 1 | 0,9715 | 31 | 0,5532 |
| 2 | 0,9540 | 32 | 0,5410 |
| 3 | 0,9374 | 33 | 0,5289 |
| 4 | 0,9215 | 34 | 0,5170 |
| 5 | 0,9061 | 35 | 0,5052 |
| 6 | 0,8911 | 36 | 0,4936 |
| 7 | 0,8764 | 37 | 0,4822 |
| 8 | 0,8619 | 38 | 0,4710 |
| 9 | 0,8477 | 39 | 0,4599 |
| 10 | 0,8336 | 40 | 0,4490 |
| 11 | 0,8196 | 41 | 0,4383 |
| 12 | 0,8056 | 42 | 0,4278 |
| 13 | 0,7918 | 43 | 0,4175 |
| 14 | 0,7780 | 44 | 0,4074 |
| 15 | 0,7642 | 45 | 0,3974 |
| 16 | 0,7505 | 46 | 0,3877 |
| 17 | 0,7368 | 47 | 0,3782 |
| 18 | 0,7232 | 48 | 0,3687 |
| 19 | 0,7096 | 49 | 0,3595 |
| 20 | 0,6961 | 50 | 0,3505 |
| 21 | 0,6826 | 51 | 0,3417 |
| 22 | 0,6692 | 52 | 0,3330 |
| 23 | 0,6559 | 53 | 0,3246 |
| 24 | 0,6427 | 54 | 0,3163 |
| 25 | 0,6295 | 55 | 0,3082 |
| 26 | 0,6165 | 56 | 0,3003 |
| 27 | 0,6036 | 57 | 0,2925 |
| 28 | 0,5908 | 58 | 0,2849 |
| 29 | 0,5781 | 59 | 0,2775 |
| 30 | 0,5656 | 60 | 0,2703 |

TABLE III
SYSTEM RELIABILITY FOR $\lambda = 0{,}05$ AND $\mu = 0{,}1$.

| t | Reliability | t | Reliability | t | Reliability | t | Reliability |
|---|---|---|---|---|---|---|---|
| 1 | 0,9501 | 21 | 0,5504 | 41 | 0,5595 | 61 | 0,3942 |
| 2 | 0,9279 | 22 | 0,5323 | 42 | 0,5504 | 62 | 0,3869 |
| 3 | 0,9067 | 23 | 0,5145 | 43 | 0,5413 | 63 | 0,3799 |
| 4 | 0,8861 | 24 | 0,4972 | 44 | 0,5323 | 64 | 0,3729 |
| 5 | 0,8659 | 25 | 0,4802 | 45 | 0,5234 | 65 | 0,3660 |
| 6 | 0,8459 | 26 | 0,4636 | 46 | 0,5145 | 66 | 0,3592 |
| 7 | 0,8259 | 27 | 0,4475 | 47 | 0,5058 | 67 | 0,3525 |
| 8 | 0,8058 | 28 | 0,4317 | 48 | 0,4972 | 68 | 0,3459 |
| 9 | 0,7857 | 29 | 0,4164 | 49 | 0,4886 | 69 | 0,3394 |
| 10 | 0,7656 | 30 | 0,4015 | 50 | 0,4802 | 70 | 0,3330 |
| 11 | 0,7454 | 31 | 0,3869 | 51 | 0,4719 | 71 | 0,3268 |
| 12 | 0,7253 | 32 | 0,3729 | 52 | 0,4636 | 72 | 0,3206 |
| 13 | 0,7052 | 33 | 0,3592 | 53 | 0,4555 | 73 | 0,3145 |
| 14 | 0,6851 | 34 | 0,3459 | 54 | 0,4475 | 74 | 0,3085 |
| 15 | 0,6652 | 35 | 0,3331 | 55 | 0,4395 | 75 | 0,3026 |
| 16 | 0,6455 | 36 | 0,3206 | 56 | 0,4317 | 76 | 0,2968 |
| 17 | 0,6259 | 37 | 0,3085 | 57 | 0,4240 | 77 | 0,2911 |
| 18 | 0,6066 | 38 | 0,2968 | 58 | 0,4164 | 78 | 0,2855 |
| 19 | 0,5876 | 39 | 0,2855 | 59 | 0,4089 | 79 | 0,2800 |
| 20 | 0,5688 | 40 | 0,2746 | 60 | 0,4015 | 80 | 0,2746 |

In Table 1-3, the time-dependent reliability of the 2-component parallel-connected technical system was calculated under the principle of repair. It was observed how the change in the reliability of the system decreased in the observed period by taking the mean of the time between failures and the time between repairs differently.

## VI. CONCLUSION

Examination of technical systems in the early literature begins with the sequence statistics of the lifetime of the components that constitute the technical system. In the following periods, auxiliary indicators that define the reliability of the technical system, such as the system signature, were introduced. System signature is an important indicator used in defining the technical system and calculating its reliability [4-7]. System signature has not been taken into consideration in repairable systems until today. This is because the system signature is formed based on the failure of components once and possible states of the technical system. Since the possible states of a technical system working under the repairable principle are not finite, it is very difficult to calculate the system signature for repairable systems. In this study, the system signature with infinite components was successfully calculated under the condition that the lifetimes of the components that constitute the technical system are exponential, and the reliability of the technical system was calculated based on the system signature (Table1-3). Research on the principle of repair or renewal in technical systems begins in the last years of the literature. In such studies, the system can be designed as a closed queue model and the reliability of the system can be calculated through the Markov model representing the system. Gurcan et al studied the repairable $k$-out of-n system under the assumption that the time elapsed between components is known and obtained two important results [8]. In the first one, the meantime to failure can be calculated. In the second, it is possible to calculate the mean working time between consecutive failures. However, unlike other studies, Gurcan et al used the common distribution of the inter-arrival failure times instead of the lifetime of the components. In this way, the difficulty encountered in determining the states of the Markov process representing the system was eliminated. In some studies based on the structure of the technical system in the literature, the results of the studies are limited because the lifetimes of the components are taken [9]. In this respect, the system signature can be used as an alternative approach in examining repairable technical systems.

## APPENDIX

Matlab code for reliability calculation which uses in the manuscript.

```
syms y
a=input('enter max value of t = ');
L=1/10; %input('enter Lambda value = '); % Lambda parameter.
M=1/5; %input('enter Mu value = '); % Mu parameter.
n=input('n = '); % value of n of Lambda
A=zeros(a,n);
for t=1:a
    A(t,1)=exp(-L*t);
end
for b=2:n
for t=1:a
  z2=0;z1=zeros(1,1);T1=0;T2=0;z=0;m=b-1;
for k=1:(b)
   T1=L^(k-1)*M^m*exp((-L)*t)*...
     int(y^(m-1)*(t-y)^(k-1)*exp((L-M)*y),'y',0,t)/...
     (factorial(k-1)*factorial(m-1));
   z1(1,1)=z1(1,1)+T1;
end
for i=1:m
  T2=(M*t)^(i-1)*exp(-M*t)/factorial(i-1);
  z2=z2+T2;
end
   A(t,b)=z1+z2;
end
  if (b==n) break
end
 end

n0=2;
for k=1:n0
 for j=1:n0
   a0(j,1)=nchoosek(n0 - 1, j-1)/(2^(n0-1));
 end
    n1=n0+1;
```

```
     for j=1:n1
         b0(j,1)=nchoosek(n1 - 1, j-1)/(2^(n1-1));
     end
         n2=n1+1;
     for j=1:n2
         c(j,1)=nchoosek(n2 - 1, j-1)/(2^(n2-1));
     end
         n3=n2+1;
     for j=1:n3
         d(j,1)=nchoosek(n3 - 1, j-1)/(2^(n3-1));
     end
         n4=n3+1;
     for j=1:n4
         e(j,1)=nchoosek(n4 - 1, j-1)/(2^(n4-1));
     end
         n5=n4+1;
     for j=1:n5
         f(j,1)=nchoosek(n5 - 1, j-1)/(2^(n5-1));
     end
         n6=n5+1;
     for j=1:n6
         g(j,1)=nchoosek(n6 - 1, j-1)/(2^(n6-1));
     end
         n7=n6+1;
     for j=1:n7
         h(j,1)=nchoosek(n7 - 1, j-1)/(2^(n7-1));
     end
         n8=n7+1;
     for j=1:n8
         l(j,1)=nchoosek(n8 - 1, j-1)/(2^(n8-1));
     end
     end
     % Probability
PT2=zeros(a,1);PT3=zeros(a,1);PT4=zeros(a,1);
PT5=zeros(a,1);PT6=zeros(a,1);PT7=zeros(a,1);
PT8=zeros(a,1);PT9=zeros(a,1);PT10=zeros(a,1);
for i=1:a
PT2(i,1)=A(i,1)*a0(1,1)+A(i,2)*a0(2,1);
PT3(i,1)=A(i,1)*b0(1,1)+A(i,2)*b0(2,1)+A(i,3)*b0(3,1);
PT4(i,1)=A(i,1)*c(1,1)+A(i,2)*c(2,1)+A(i,3)*c(3,1)+…
    A(i,4)*c(4,1);
PT5(i,1)=A(i,1)*d(1,1)+A(i,2)*d(2,1)+A(i,3)*d(3,1)+…
    A(i,4)*d(4,1)+A(i,5)*d(5,1);
PT6(i,1)=A(i,1)*e(1,1)+A(i,2)*e(2,1)+A(i,3)*e(3,1)+…
    A(i,4)*e(4,1)+A(i,5)*e(5,1)+A(i,6)*e(6,1);
PT7(i,1)=A(i,1)*f(1,1)+A(i,2)*f(2,1)+A(i,3)*f(3,1)+…
    A(i,4)*f(4,1)+A(i,5)*f(5,1)+A(i,6)*f(6,1)+A(i,7)*f(7,1);
PT8(i,1)=A(i,1)*g(1,1)+A(i,2)*g(2,1)+A(i,3)*g(3,1)+…
    A(i,4)*g(4,1)+A(i,5)*g(5,1)+A(i,6)*g(6,1)+…
    A(i,7)*g(7,1)+A(i,8)*g(8,1);
PT9(i,1)=A(i,1)*h(1,1)+A(i,2)*h(2,1)+A(i,3)*h(3,1)+…
    A(i,4)*h(4,1)+A(i,5)*h(5,1)+A(i,6)*h(6,1)+…
    A(i,7)*h(7,1)+A(i,8)*h(8,1)+A(i,9)*h(9,1);
PT10(i,1)=A(i,1)*l(1,1)+A(i,2)*l(2,1)+A(i,3)*l(3,1)+…
    A(i,4)*l(4,1)+A(i,5)*l(5,1)+A(i,6)*l(6,1)+A(i,7)*l(7,1)+…
    A(i,8)*l(8,1)+A(i,9)*l(9,1)+A(i,10)*l(10,1);
end
P=[PT2 PT3 PT4 PT5 PT6 PT7 PT8 PT9 PT10];
    % Reliability
R1=zeros(a,1);R=zeros(a,1);
ksi=M/(L+M);
S=[1-ksi (1-ksi)*ksi (1-ksi)*ksi^2 (1-ksi)*ksi^3…
    (1-ksi)*ksi^4 (1-ksi)*ksi^5 (1-ksi)*ksi^6…
    (1-ksi)*ksi^7 (1-ksi)*ksi^8];
for i=1:a
for j=2:n
    R1(i,1)=S(j-1)*P(i,j-1);
    R(i,1)=R(i,1)+R1(i,1);
end
end
```